\documentclass[11pt]{article}
\usepackage{amsmath,amsfonts,amssymb}
\newtheorem{theorem}{Theorem}
\newtheorem{lemma}{Lemma}

\newcommand{\R}{\mathbb{R}}

\newcommand{\N}{{\mathbb N}}

\newcommand{\E}{{\mathsf E}}

\newcommand{\Pp}{\mathsf{P}}

\newcommand{\ONE}{\mathbf{1}}

\title{Solutions of Semilinear Wave Equation via Stochastic Cascades }
\author{Yuri Bakhtin\thanks{
School of Mathematics, Georgia Institute of Technology, Atlanta, GA, 30067-0160, USA, 
  email:{ bakhtin@math.gatech.edu}, phone: +1(404)894-9235, fax: +1(404)894-4409 (corresponding author), supported by
an NSF CAREER grant.
} \and Carl Mueller\thanks{Department of Mathematics, University of Rochester, 
Rochester, NY  14627, web page: http://www.math.rochester.edu/people/faculty/cmlr/, supported by an NSF grant.}}
\date{}

\begin{document}
\maketitle
\begin{abstract}{We introduce a probabilistic representation for solutions of quasilinear wave equation with analytic
nonlinearities.
We use stochastic
cascades to prove existence and uniqueness of the solution.}

{\it Keywords:  semilinear wave equation, probabilistic representation, stochastic cascade}
\end{abstract}

We consider the following nonlinear wave equation on the real line $\R$:
\begin{equation}
\Box u(x,t)=F(x,t, u(x,t)),\quad(x,t)\in\R\times \R_+. 
\label{eq:nonlinear_wave}
\end{equation}
Here \[\Box u(x,t)=u_{tt}(x,t)-u_{xx}(x,t),\] and  
$F$ is a given function.

For $T>0$, we say that $u:\R\times [0,T)$ is a classical solution of the Cauchy problem 
of~\eqref{eq:nonlinear_wave} with initial conditions
\begin{align}
 u(x,0)&=\phi(x),\quad x\in\R,\label{eq:initial_condition_0}\\
 u_t(x,0)&=\psi(x),\quad x\in \R,\label{eq:initial_condition_1}
\end{align}
if $u\in C^2(\R\times(0,T))\cap C^1(\R\times(0,T))$.

The goal of this note is to develop a stochastic cascade approach to constructing
 solutions of the Cauchy problem
\eqref{eq:nonlinear_wave}--\eqref{eq:initial_condition_1}. It is similar to the
construction of solutions for the Navier--Stokes system suggested in \cite{LeJan-Sznitman} and for
the KPP equation in~\cite{McKean}.
Although this approach is essentially equivalent to a Picard-type iteration scheme, it provides an
interesting point of view.   

Probabilists have long desired a probabilistic representation of the wave 
equation, but there are only a few papers on the topic.  In \cite{Kac1,Kac2}, 
Kac discovered a probabilistic representation for the telegrapher's equation, 
which is a wave equation with a lower-order time derivative.  More recently, 
Dalang, Tribe, and the second author developed a multiplicative version of the 
Feynman-Kac formula which applies to the wave equation, among others, see 
\cite{dm09,dmt08}. 

We begin with the classical d'Alembert representation of solutions for the linear wave equation.
If $F\equiv0$, i.e., the problem \eqref{eq:nonlinear_wave}--\eqref{eq:initial_condition_1} 
is a homogeneous linear problem, and existence and uniqueness
hold under unrestrictive assumptions on the regularity of initial conditions $\phi$ and $\psi$. Fixing
$\phi$ and $\psi$, one can write the solution of the homogeneous Cauchy problem as
\begin{equation}
 \label{eq:homo}
  v(x,t)=\frac{1}{2}\int_{x-t}^{x+t}\psi(y)dy+\frac{1}{2}(\phi(x+t)+\phi(x-t)),
\end{equation}

If $F(x,t,u)=f(x,t)$ is a sufficiently smooth function that does not depend on 
$u$, then we have an inhomogeneous wave equation with external source
$f$, and the d'Alembert formula holds:
\begin{equation}
\label{eq:DAlembert}
 u(x,t)=v(x,t)+\frac{1}{2}\int_{\Delta(x,t)} f(x,t)\,dx\,dt,
\end{equation}
where 
\[
\Delta(x,t)=\{(y,s):\ 0\le s\le t,\ |y-x|\le t-s\} 
\]
 is the light cone of the past associated with the space-time point $(x,t)$. 
In fact, in our one-dimensional situation, $\Delta(x,t)$ is just a triangle. 

Formula \eqref{eq:DAlembert} allows us to define a mild solution of equation~\eqref{eq:nonlinear_wave}
on a time interval $[0,T)$ as a measurable function 
$u:\R\times[0,T)\to\R$ such that for all $(x,t)\in \R\times[0,T)$,
\begin{equation}
\label{eq:DAlembert_nonlin}
 u(x,t)=v(x,t)+\frac{1}{2}\int_{\Delta(x,t)} F(x,t,u(x,t))\,dx\,dt.
\end{equation}

From now on we shall assume for simplicity that $F(x,t,u)=F(u)$ does not depend on $(x,t)$, although our construction
can be also applied with appropriate modifications in the general case.
The next assumption is crucial for our construction, though: we require analyticity of $F$, i.e., 
we assume that for all $u$, $F(u)$ can be represented as a convergent power series in $u$:
\[
 F(u)=\sum_{k=0}^\infty a_ku^k.
\]
In particular, we can deal with power-type nonlinearities like $F(x,t,u)\equiv u^2$.

Let us fix $(x,t)\in\R\times\R_+$ and try to rewrite the d'Alembert formula in the language of random variables. 
To that end, let us introduce a random point $(\xi,\tau)$ uniformly distributed in $\Delta(x,t)$. Since
the area (Lebesgue measure) of $\Delta(x,t)$ equals $t^2$, it means that the random point 
$(\xi,\tau)$ has density vanishing outside of $\Delta(x,t)$ and identically equal to $t^{-2}$ inside 
$\Delta(t,x)$.
Therefore, formula~\eqref{eq:DAlembert_nonlin} can be rewritten as
\begin{equation}
 \label{eq:probabilistic_DAlembert_nonlin}
 u(x,t)=v(x,t)+\frac{t^2}{2}\E F(u(\xi,\tau)).
\end{equation}

Next step is to consider  a sequence
of numbers $(p_k)_{k=0}^\infty$ with the following properties:
\begin{enumerate}
 \item[(i)] $p$ is a probability distribution: $\sum_{k=0}^\infty p_k=1$, $p_k\ge0$ for all $k\ge 0$;
 \item[(ii)] $p_k>0$ for every $k\ge 0$ with $a_k\ne 0$;
 \item[(iii)] $\sum_{k=1}^\infty kp_k\le 1$. 
\end{enumerate}
Also let us define 
\[
w(x,t)=\frac{v(x,t)}{p_{0}},\quad (x,t)\in\R\times\R_+. 
\]
and, for each $k\ge 0$,
\[
 b_k=\begin{cases}\frac{a_k}{p_k},&  p_k\ne 0, \\0,& p_k=0,\end{cases}
\]

Let us introduce a random variable $\kappa$ distributed according to $p$ and independent of $(\xi,\tau)$. Then
\eqref{eq:probabilistic_DAlembert_nonlin} immediately implies the following lemma.
\begin{lemma}\label{lm:iteration_step} If $u$ is a solution of~\eqref{eq:nonlinear_wave} on $[0,T)$, then for any
$(x,t)\in\R\times\R_+$,
\begin{equation}
 \label{eq:probabilistic_DAlembert_nonlin1}
 u(x,t)=\E
\left[\left(w(x,t)+\frac{t^2b_0}{2}\right)\ONE_{\{\kappa=0\}}+\frac{t^2}{2}b_{\kappa}u^{\kappa}(\xi,\tau)\ONE_{\{
\kappa\ge 1\}}\right].
\end{equation}
\end{lemma}

Next step in our construction is to iterate Lemma~\ref{lm:iteration_step}.  
Namely, for any $k\ge 0$, on the event $\{\kappa=k\}$ we may compute 
$u^{k}(\xi,\tau)$ by the same procedure. The role of $(x,t)$ is played by 
$(\xi,\tau)$, and, conditioned on $\{(\xi,\tau)=(x',t')\}$, to compute the 
product of $k$ copies of $u(x',t')$ we may consider $k$ independent random 
variables $(\xi_{i},\tau_{i},\kappa_i)_{i=1}^k$, so that random points 
$(\xi_i,\tau_i)$ are uniformly distributed in $\Delta_{x',t'}$, and random 
variables $\kappa_i$ are distributed according to distribution $p$. Given that 
collection of random variables, for each $i=1,\ldots,k$, we can apply 
Lemma~\ref{lm:iteration_step}. Notice that on the event $\{\kappa=0\}$, the 
random variable under the expectation sign in~\eqref{eq:probabilistic_DAlembert_nonlin1}  is 
a constant equaling $w(x,t)+\frac{t^2b_0}{2}$, so that we do not have to 
consider any new random variables to compute it. 

It is clear that iterating this procedure we obtain a stochastic cascade, 
i.e., a branching process with each particle assigned a location in space and 
time.  To make this idea precise, let us introduce more notation. We shall need a 
probability space rich enough to support these random structures involving 
random family trees of the participating particles and their random locations.

We begin with an encoding of vertices of finite rooted ordered trees. Each vertex
$v$ in the $n$-th generation of the tree can be identified with a sequence $(a_1,\ldots, a_n)$, where
$a_i\in\{0,1,\ldots\}$ for all $i=1,\ldots,n$. The parent of $(a_1,\ldots, a_n)$ is $(a_1,\ldots, a_{n-1})$. It is
convenient to identify the root of the tree with an empty sequence denoted by $\emptyset$ which is consistent with the
above encoding of the parent-child relation.

Let us fix $n\in\N$ and let probability measure on rooted ordered trees with at most $n$ generations be given by the
classical Galton--Watson distribution on trees with branching distribution $(p_k)_{k=0}^{\infty}$, i.e., each vertex
$v$ of the tree has a random number $\kappa_v$ of children $(v,1),\ldots, (v,\kappa)$, where $(v,i)$ means a sequence
obtained from $v$ by appending (concatenating) $i$ on the right. The random variables $k_v$ are jointly independent.

Each vertex $v$
in thus generated random tree gets a random space-time label $(\kappa_v,\tau_v)$ according to the following
rule. First, we set $(\kappa_\emptyset,\tau_\emptyset)=(x,t)$. Then, we can iteratively apply the following: for any
vertex $v=(a_1,\ldots,a_m)$, $m< n$, conditioned on $\kappa_{v},\tau_{v}=(y,s)$ and on $\kappa_{v}=k$, the labels
$(\kappa_{(v,i)},\tau_{(v,i)})_{i=1}^k$ are
i.i.d. uniform random points in $\Delta(y,s)$, independent of all previously constructed random elements in the procedure.

Now we shall recursively define a functional $\Pi$ on subtrees of individual vertices starting with the leaves of the
tree. There are two types of the leaves. Leaves of type 1 are vertices of generation $n$. For any leaf $v$ of type~1,
we set $\Pi(v)=u(\xi_v,\tau_v)$. A leave of type 2 is a vertex $v$ that did not produce any children, i.e.,
$\kappa_v=0$. For these vertices we set 
\begin{equation}
\label{eq:leaves}
\Pi(v)=w(\xi_v,\tau_v)+\frac{\tau_v^2b_0}{2}. 
\end{equation}
As soon as the values of $\Pi(v)$ have been assigned for all leaves of the realization of the tree, we
may start assigning values to other vertices iteratively. For any vertex $v$ 
of the tree such that $\Pi(v,i)$ has already been assigned for all $i=1,\ldots,\kappa_v$, we define
\begin{equation}
\label{eq:non-leaves}
 \Pi(v)= \frac{\tau_v^2}{2}b_{\kappa_v }\prod_{i=1}^{\kappa_v}u(\xi_{(v,i)},\tau_{(v,i)}).
\end{equation}
Proceeding iteratively, we eventually will assign some value $\Pi(\emptyset)$ to the root of the tree. This value is a
random variable incorporating information from all the other vertices of the random tree as well as their
space-time labels. We denote the resulting random variable by $\Pi_n(x,t,u(\cdot))$ to stress that the tree of
at maximum $n$ generations was initiated at $(x,t)$ and the solution $u$ was used to evaluate $\Pi$ at the tree's
leaves.
\begin{theorem}
 \label{th:finite_iteration_scheme}
Suppose $u$ is a solution of~\eqref{eq:nonlinear_wave} on $[0,T]$. Then for any $n$ and any $(x,t)\in\R\times (0,T)$,
\[
 u(x,t)=\E \Pi_n(x,t,u(\cdot)).
\]
\end{theorem}

An obvious next step is to take $n$ to infinity. Since the requirement (iii) on the branching distribution $p$ means
that the branching process is critical or subcritical, the realizations of the random trees almost surely have finitely
many vertices. In particular, with probability 1,
\begin{equation}
 \label{eq:convergence_of_rv}
 \lim_{n\to\infty} \Pi_n(x,t,u(\cdot)) = \Pi_\infty(x,t),
\end{equation}
where the random variable $\Pi_\infty(x,t)$ is constructed from the realization of the stochastic cascade in exactly the
same way as $\Pi_n(x,t,u(\cdot))$ for finite $n$ except that there are no leaves of type~1.

If we show that 
\begin{equation}
 \label{eq:convergence_of_expectations}
 \lim_{n\to\infty} \E \Pi_n(x,t,u(\cdot)) = \E \Pi_\infty(x,t),
\end{equation}
then we will be able to conclude that 
\begin{equation}
\label{eq:solution_representation}
u(x,t)=\E \Pi_\infty(x,t).
\end{equation}
In that case, since 
$\Pi_\infty(x,t)$ is a (random) functional of $w$, the (modified) external source only and does not involve $u$, we can
claim that the solution is unique and it is given by formula~\eqref{eq:solution_representation}. It is also easy to see
from~\eqref{eq:probabilistic_DAlembert_nonlin1} that if~\eqref{eq:convergence_of_expectations} holds then $u(x,t)$ given
by formula~\eqref{eq:solution_representation} is a solution.

So, condition~\eqref{eq:convergence_of_expectations} implies the existence and uniqueness of solution and its
stochastic representation.

We notice that $\Pp\{\Pi_\infty\ne \Pi_n\}\to 0$ as $n\to\infty$. Therefore, to ensure convergence in
~\eqref{eq:convergence_of_expectations} it is sufficient to check that
$\E \Pi_\infty(x,t)$ is well-defined.

Next, we notice that, according to~\eqref{eq:leaves} and~\eqref{eq:non-leaves} $\Pi(x,t)$ is a product of many
factors of the form $w(\xi_v,\tau_v)+\frac{\tau_v^2b_0}{2}$ or $\frac{\tau_v^2b_{\kappa_v}}{2}$. If we require
that all these factors are bounded by $1$, then $\Pi_\infty(x,t)$ product is bounded by 1, and its expectation is also
bounded. This leads us to the following theorem.

Let $b_*=\sup_{k\in\N}|b_k|$. Let
\[
 T^*=\sup\left\{t\le \sqrt{\frac{2}{b^*}}:\ \sup_{x\in\R} |w(x,t)|+\frac{t^2b_0}{2}\le 1\right\}.
\]

\begin{theorem}\label{th:boundedness_conditions} If $T^*>0$,
then there is a unique solution of \eqref{eq:nonlinear_wave} on $[0,T^*)$. It is given by
\eqref{eq:solution_representation}.
\end{theorem}

The crude requirement of boundedness by 1 and the resulting condition of the above theorem can certainly be improved
for some specific cases. We do not explore
this issue further since we do not expect our method to produce sharp conditions for the existence of the solution.

\bibliographystyle{alpha}

\end{document}